\documentclass{amsart}
\newtheorem{thm}{Theorem}
\newtheorem{cor}[thm]{Corollary}
\newcommand{\G}{\Gamma} 
\newcommand{\U}{{\mathcal U}} 
\newcommand{\e}{\varepsilon} 
\newcommand{\N}{{\mathbb N}} 
\newcommand{\B}{{\mathbb B}} 
\newcommand{\M}{{\mathbb M}} 
\newcommand{\R}{{\mathbb R}} 
\newcommand{\F}{{\mathbb F}} 
\newcommand{\hh}{{\mathcal H}} 
\newcommand{\id}{\mathrm{id}} 
\begin{document}
\title{There is no separable universal $\mathrm{II}_1$-factor}
\author{Narutaka OZAWA}
\address{Department of Mathematical Science, 
University of Tokyo, 153-8914, Japan}
%\curraddr{}
\email{narutaka@ms.u-tokyo.ac.jp}
\thanks{Partially supported by JSPS Postdoctoral 
Fellowships for Research Abroad.}
\date{November 1, 2002.}
\begin{abstract}
Gromov constructed uncountably many pairwise non-isomorphic 
discrete groups with Kazhdan's property $\mathrm{(T)}$. 
We will show that no separable $\mathrm{II}_1$-factor can contain 
all these groups in its unitary group. 
In particular, no separable $\mathrm{II}_1$-factor can contain 
all separable $\mathrm{II}_1$-factors in it. 
We also show that the full group $C^*$-algebras of 
some of these groups fail the lifting property. 
\end{abstract}
\maketitle
\section{Results}
We recall that a discrete group $\G$ 
is said to have Kazhdan's property $\mathrm{(T)}$ 
if the trivial representation is isolated in the dual $\hat{\G}$ of $\G$, 
equipped with the Fell topology. 
This is equivalent to that there are 
a finite subset $E$ of generators in $\G$ and 
a decreasing function $f\colon \R_+\to\R_+$ 
with $\lim_{\e\to 0}f(\e)=0$ such that the following is true: 
if $\pi$ is a unitary representation of $\G$ on a Hilbert space $\hh$ 
and $\xi\in\hh$ is a unit vector with 
$\e=\max_{s\in E}\|\pi(s)\xi-\xi\|$, 
then there is a vector $\eta\in\hh$ with 
$\|\xi-\eta\|<f(\e)$ (in particular $\eta\neq0$ when $\e$ is small enough)
such that $\pi(s)\eta=\eta$ for all $s\in\G$. 
We refer the reader to \cite{hv} and \cite{valette} for 
the information of Kazhdan's property $\mathrm{(T)}$.
We recall that a discrete group $\G$ is said to be quasifinite 
if all its proper subgroups are finite, 
and is said to be infinite conjugacy classes (abbreviated to ICC) 
if all nontrivial conjugacy classes in $\G$ are infinite. 
We note that a discrete group $\G$ is ICC 
if and only if its group von Neumann algebra $L\G$ is a factor. 
We also observe that a group which is quasifinite and ICC has to be simple. 

Gromov (Corollary 5.5.E in \cite{gromov}) claimed 
that any torsion-free non-cyclic hyperbolic group has 
a quotient group all of whose proper subgroups are cyclic 
of prescribed orders (cf.\ Theorem 3.4 in \cite{valette}). 
This claim was partly confirmed by Olshanskii 
(Corollary 4 in \cite{olshanskii}). 
Actually, what Olshanskii proved there is that 
any torsion-free non-cyclic hyperbolic group has 
a nontrivial quasifinite quotient group. 
We observe that Olshanskii's argument gives us the following. 
\begin{thm}[Gromov-Olshanskii]\label{hyp}
Any torsion-free non-cyclic hyperbolic group has 
uncountably many pairwise non-isomorphic quotient groups 
all of which are quasifinite and ICC. 
In particular, there is a discrete group $\G$ 
with Kazhdan's property~$\mathrm{(T)}$ which has 
uncountably many pairwise non-isomorphic quotient groups 
$\{\G_\alpha\}_{\alpha\in I}$ all of which are simple and ICC. 
\end{thm}
Connes conjectured that a discrete group $\Delta$ with 
Kazhdan's property $\mathrm{(T)}$ and the ICC property 
is uniquely determined by its group von Neumann algebra $L\Delta$. 
The following theorem and its corollary, 
which was suggested by S. Popa, 
confirm Connes' conjecture for $\{\G_\alpha\}_{\alpha\in I}$ 
``modulo countable sets'' and solve Problem 4.4.29 in \cite{sakai}, 
Conjecture 4.5.5 in \cite{popa} and Problem III.45 in \cite{delaharpe}. 
See also Theorem 1 in \cite{popa2} and its remarks. 
\begin{thm}\label{vna}
Let $\G$ and $\{\G_\alpha\}_{\alpha\in I}$ be as in Theorem \ref{hyp} 
and let $M$ be a separable $\mathrm{II}_1$-factor.
Then, the set 
$$\{ \alpha\in I : \mbox{the unitary group $\U(M)$ of $M$ 
contains a subgroup isomorphic to $\G_\alpha$}\}$$
is at most countable. 
\end{thm}
Recall that two $\mathrm{II}_1$-factors $M$ and $N$ 
are said to be stably equivalent if there are $n\in\N$ and 
a projection $p\in\M_n(M)$ such that $p\M_n(M)p$ is isomorphic to $N$. 
\begin{cor}\label{corv}
Let $\G$ and $\{\G_\alpha\}_{\alpha\in I}$ be as in Theorem \ref{hyp} 
and let $M$ be a separable $\mathrm{II}_1$-factor.
Then, the set 
$$\{ \alpha\in I : \mbox{$M$ contains a subfactor 
which is stably equivalent to $L\G_\alpha$}\}$$
is at most countable. 
\end{cor}
In connection with Connes' embedding problem \cite{connes}, 
it would be interesting to know whether all (or at least one of) 
$\G_\alpha$'s are embeddable into the unitary group $\U(R^\omega)$ 
of the ultrapower $R^\omega$ of hyperfinite $\mathrm{II}_1$-factors. 
Since each $\G_\alpha$ arises as a limit of hyperbolic groups, 
we observe that if all hyperbolic groups are embeddable into 
$\U(R^\omega)$, then so is $\G_\alpha$. 
We remark that whether all hyperbolic groups are residually finite 
(and thus embeddable into $\U(R^\omega)$) is one of the major 
open problems in geometric group theory. 

Let us consider the category of unital $C^*$-algebras 
and unital completely positive maps. 
A $C^*$-algebra $A$ is said to be complementary universal 
for a class $\mathcal{C}$ of $C^*$-algebras 
if for every member $B$ of $\mathcal{C}$, 
there are unital completely positive maps 
$\psi\colon B\to A$ and $\varphi\colon A\to B$ such that 
$\varphi\psi=\id_B$. 
It follows from Kirchberg's theorem \cite{kirchberg}
that any separable $C^*$-algebra of not type I 
is complementary universal for the class of separable 
nuclear $C^*$-algebras. 
The full group $C^*$-algebra $C^*\F_\infty$ of 
the free group $\F_\infty$ on countably many generators 
is complementary universal for the class of separable 
$C^*$-algebras with the lifting property (abbreviated to LP). 
See \cite{klp} for the information of the LP. 
It is not known whether there exists a separable 
complementary universal $C^*$-algebra for the class of 
separable exact $C^*$-algebras. 
\begin{thm}\label{cst}
Let $\G$ and $\{\G_\alpha\}_{\alpha\in I}$ be as in Theorem \ref{hyp} 
and let $\mathcal{C}=\{C^*\G_\alpha : \alpha\in I\}$ or 
$\mathcal{C}=\{C^*_{\mathrm{red}}\G_\alpha : \alpha\in I\}$.
Then, there is no separable complementary universal $C^*$-algebra 
for $\mathcal{C}$. 
\end{thm}
 From the above discussion, we immediately obtain 
the following corollary. 
\begin{cor}
The full group $C^*$-algebra $C^*\G_\alpha$ of $\G_\alpha$ 
fails the LP for some $\alpha\in I$. 
\end{cor}
\section*{Acknowledgment}
The author thanks Professor S. Popa for useful comments. 
This research was carried out while the author was visiting 
the University of California Berkeley 
under the support of the Japanese Society for the Promotion of 
Science Postdoctoral Fellowships for Research Abroad. 
He gratefully acknowledges the kind hospitality from UCB. 
\section{Proofs}
\begin{proof}[Proof of Theorem \ref{hyp}]
Since there exists a torsion-free non-cyclic hyperbolic group 
with Kazhdan's property $\mathrm{(T)}$ 
(e.g., a co-compact lattice in $Sp(n,1)$ or in $F_{4(-20)}$), 
the second part is a straight consequence of the first. 
We just indicate how to modify 
the proof of Corollary 3 in \cite{olshanskii} 
to obtain the first part of our Theorem \ref{hyp}. 
So, we stick by the notations used in \cite{olshanskii}. 
Let $G=\{g_1,g_2,\ldots\}$ be a torsion-free non-cyclic hyperbolic group. 
It follows that $G$ is ICC since 
the set $E(G)$ of all $x\in G$ whose conjugacy class is finite 
is a finite subgroup in $G$ (cf.\ Proposition 1 in \cite{olshanskii}). 
Recall that the quasifinite quotient group $G'$ was 
the inductive limit of a sequence of epimorphisms 
$G=G_0\to G_1\to G_2\to\cdots$. 
By the construction, every $E(G_i)$ is trivial, 
or equivalently every $G_i$ is ICC. 
Hence, manipulating the construction, for every $i$ and $j\le i$, 
we can carry at least $i$ mutually distinct elements 
from the conjugacy class of $g_j$ in $G_i$ injectively
into $G'$ unless $g_j=1$ in $G_i$. 
This ensures the ICC property of $G'$. 
In the construction, there has to be infinitely many $i$'s 
such that $g_i$ is torsion-free in $G_{2i-2}$. 
For such $i$, we may choose arbitrarily large number for 
the order of $g_i$ in $G_{2i-1}$ which will be equal to that in $G'$. 
Combined with a diagonal argument, this implies 
that there are uncountably many normal subgroups in $G$ all of 
whose corresponding quotient groups are quasifinite and ICC. 
Theorem~\ref{hyp} now follows from this result 
and Lemma III.42 in \cite{delaharpe}. 
\end{proof}
\begin{proof}[Proof of Theorem \ref{vna}]
To prove the theorem by reductio ad absurdum, 
suppose that 
$$I_0=\{ \alpha\in I : \mbox{$\U(M)$ contains 
a subgroup isomorphic to $\G_\alpha$}\}$$
is uncountable. 
For each $\alpha\in I_0$, let $u_\alpha\colon\G\to\U(M)$ be 
a non-trivial homomorphism which factors through $\G_\alpha$. 
We fix a standard representation of $M$ on $\hh$ with 
a unit cyclic separating trace vector $\xi$ in $\hh$.  
It follows that there are $\delta>0$ and 
an uncountable subset $I_1$ of $I_0$ such that 
$\max_{s\in\G}\|u_\alpha(s)\xi-\xi\|>\delta$ for all $\alpha\in I_1$. 

Let a finite subset $E$ of generators in $\G$ 
and a function $f$ be as in the above definition of 
Kazhdan's property $\mathrm{(T)}$. 
Take $\e>0$ small enough so that $2f(\e)<\delta$. 
Since $\hh$ is separable and $I_1$ is uncountable, 
there are distinct $\alpha$ and $\beta$ in $I_1$ such that 
$\max_{s\in E}\|u_\alpha(s)\xi-u_\beta(s)\xi\|<\e$. 
We consider the unitary representation 
$\pi\colon\G\ni s\mapsto u_\alpha(s)Ju_\beta(s)J\in\B(\hh)$, 
where $J$ is the canonical conjugation on $\hh$ associated 
with $M$ and $\xi$.
Then, we have $\max_{s\in E}\|\pi(s)\xi-\xi\|<\e$. 
It follows from Kazhdan's property $\mathrm{(T)}$ of $\G$ 
that there is a vector $\eta\in\hh$ with $\|\xi-\eta\|<f(\e)$ 
such that $\pi(s)\eta=\eta$ for all $s\in\G$. 
Let $\Delta=\{ s\in\G : u_\alpha(s)\eta=\eta \}$. 
It is easy to see that $\Delta$ is a subgroup of $\G$ 
and that $\Delta$ contains 
the normal subgroups $\ker u_\alpha$ and $\ker u_\beta$.
Since $\G_\alpha$ and $\G_\beta$ are simple and 
$\ker u_\alpha$ and $\ker u_\beta$ are distinct, 
we actually have $\Delta=\G$. 
It follows that $\max_{s\in\G}\| u_\alpha(s)\xi-\xi\|<2f(\e)<\delta$, 
which is absurd. 
\end{proof}
\begin{proof}[Proof of Corollary \ref{corv}]
It is not difficult to see that 
if $L\G_\alpha$ is isomorphic to a (not necessarily unital) 
subfactor of $M$, 
then $\G_\alpha$ is isomorphic to a subgroup of $\U(M)$. 
Therefore, it follows from Theorem \ref{vna} that 
$$\{ \alpha\in I : \mbox{$\M_n(M)$ contains a (not necessarily unital) 
subfactor isomorphic to $L\G_\alpha$}\}$$
is at most countable for every $n\in\N$, and the conclusion follows. 
\end{proof}
\begin{proof}[Proof of Theorem \ref{cst}]
We only deal with the case 
where $\mathcal{C}=\{ C^*\G_\alpha : \alpha\in I\}$. 
To prove the theorem by reductio ad absurdum, 
suppose that there is a separable $C^*$-algebra $A$ 
which is complementary universal for $\mathcal{C}$. 
We fix unital completely positive maps 
$\psi_\alpha\colon C^*\G_\alpha\to A$ and 
$\varphi_\alpha\colon A\to C^*\G_\alpha$ such that 
$\varphi_\alpha\psi_\alpha=\id_{C^*\G_\alpha}$.
Let $E\subset\G$ be a finite set of generators of $\G$ 
containing $1$ and let $u_\alpha(s)$ be 
the unitary element in $C^*\G_\alpha$ corresponding to $s\in\G$. 
Let $\e>0$ be arbitrary. 
Since $A$ is separable and $I$ is uncountable, 
there are distinct $\alpha$ and $\beta$ in $I$ such that 
$\max_{s\in E}\|\psi_\alpha(u_\alpha(s))-\psi_\beta(u_\beta(s))\|<\e$.
It follows that denoting the left regular representation 
of $\G_\alpha$ by $\lambda_\alpha$, we have 
\begin{align*}
\|\frac{1}{|E|} & 
 \sum_{s\in E}\lambda_\alpha(u_\alpha(s))\otimes u_\beta(s)
 \|_{C^*_{\mathrm{red}}\G_\alpha\otimes_{\max}C^*\G_\beta}\\
&\geq\| \frac{1}{|E|}\sum_{s\in E}
 \lambda_\alpha(u_\alpha(s))\otimes\varphi_\alpha\psi_\beta(u_\beta(s))
 \|_{C^*_{\mathrm{red}}\G_\alpha\otimes_{\max}C^*\G_\alpha}\\
&\geq\|\frac{1}{|E|}\sum_{s\in E}
 \lambda_\alpha(u_\alpha(s))\otimes u_\alpha(s)
 \|_{C^*_{\mathrm{red}}\G_\alpha\otimes_{\max}C^*\G_\alpha}-\e\\
&=1-\e.
\end{align*}
Since $\G$ has Kazhdan's property $\mathrm{(T)}$, if we choose 
$\e>0$ sufficiently small, this implies that 
the trivial representation of $\G$ is weakly contained in 
$C^*_{\mathrm{red}}\G_\alpha\otimes_{\max}C^*\G_\beta$ 
(cf.\ Proposition 4.9 in \cite{valette}). 
Reasoning in the same way as the proof of Theorem~\ref{vna}, 
one can show that the trivial representation is weakly 
contained in $C^*_{\mathrm{red}}\G_\alpha$. 
This is absurd. 
\end{proof}

\end{document}